\documentclass[12pt]{article}

\usepackage{latexsym,epsfig}
\usepackage{amssymb,amsmath,amsthm}

\newcommand{\la}{\lambda}

\begin{document}
\title{Existence
of sign-changing solutions for the nonlinear $p$-Laplacian boundary value problem}
\author
{Wei-Cheng Lian$^\dag$, Wei-Chuan Wang$^\ddag$ and Y.H. Cheng$^*$}
\maketitle
 \begin{abstract}
We study the nonlinear one-dimensional $p$-Laplacian equation $$
-(y'^{(p-1)})'+(p-1)q(x)y^{(p-1)}=(p-1)w(x)f(y) \mbox{ on }
(0,1),$$ with linear separated boundary conditions. We give
sufficient conditions for the existence of solutions with
prescribed nodal properties concerning the behavior of
$f(s)/s^{(p-1)}$ when $s$ are at infinity and zero. These results are more
general and complementary for previous known ones for the case
$p=2$ and $q$ is nonnegative.
\end{abstract}
 \vskip0.5in
AMS Subject Classification (2000) : 34A55, 34B24.

\noindent Keywords: Existence, zero, solutions, nonlinear,
$p$-Laplacian, boundary value problem.

 \footnote{ $^\dag$Department of Information
Management, National
 Kaohsiung Marine Univeristy, Kaohsiung 811, Taiwan, R.O.C. Email:
 wclian@mail.nkmu.edu.tw}
\footnote{ $^\ddag$Department of Applied Mathematics, National Sun Yat-sen
University, Kaohsiung, Taiwan 804, R.O.C. Email:
wangwc@math.nsysu.edu.tw}
\footnote{$^*$Department of Mathematics, National Tsing Hua University, Hsinchu, Taiwan 300, R.O.C.
Email:jengyh@math.nsysu.edu.tw}

\def\theequation{\arabic{section}.\arabic{equation}}
\section{Introduction}
\setcounter{equation}{0} \hskip0.25in Consider the following
nonlinear one-dimensional $p$-Laplacian equation,
\begin{equation}
\label{eq1.1}
-(y'^{(p-1)})'+(p-1)q(x)y^{(p-1)}=(p-1)w(x)f(y) \mbox{ on } (0,1),
\end{equation}
with the separated boundary conditions
\begin{eqnarray}\left\{\begin{array}{l}\label{eq1.2}
S_p'(\alpha)y(0)-S_p(\alpha)y'(0)=0,~~\alpha\in[0,\pi_p),\\
S_p'(\beta)y(1)-S_p(\beta)y'(1)=0,~~\beta\in(0,\pi_p],
\end{array}\right.\
\end{eqnarray}
where $p>1, \displaystyle y^{(p-1)}=\left|y\right|^{p-1}
\mbox{sgn}~y=\left|y\right|^{p-2}y$.

Denote by $w=w(x)=S_p (x)$ the inverse function of the integral
$$x=\int_{0}^{w}\frac{dt}{(1-t^{p})^{\frac{1}{p}}}\ ,\
\mbox{for}\ 0\leq w\leq 1\ .$$
In 1979, Elbert \cite{e79} discussed the analogies between $S_p
(x)$ and $\sin x$. He showed that $S_p (x)$ is  the solution of
\begin{equation}
\label{eq1.25}
\left\{\begin{array}{l} (y'^{(p-1)})'=-(p-1)y^{(p-1)}\ ,\\ y(0)=0\ ,\ y'(0)=1\ , \end{array}\right.
\end{equation}
and $\pi_{p}\equiv
\int_{0}^{1}\frac{2dt}{(1-t^{p})^{\frac{1}{p}}}=\frac{2\pi/p}{
\sin(\pi/p)}$ is the
first zero of $S_p  (x)$. Furthermore, defining
$$S_p (x)=\left\{\begin{array}{ll}S_{p}(\pi_{p}-x)\ ,& \mbox{if}\ \frac{\pi_{p}}{2}\leq x\leq \pi_{p}\ ,\\
-S_{p}(x-\pi_{p})\ ,& \mbox{if}\ \pi_{p}\leq x\leq 2\pi_{p}\ ,\\
S_{p}(x-2n\pi_{p})\ ,& \mbox{for}\ n=\pm 1, \pm 2, \cdots\ ,\\
 \end{array} \right.$$
he obtained a sine-like function and the function $S_{p}$ is so
called the generalized sine function.

The application of the most original authors cited nowadays is the
highly viscid fluid flow (cf. Ladyzhenskaya
\cite{la}, Lions \cite{li})
\begin{eqnarray*}(BVP_1)\left\{\begin{array}{l}
 -div \big((\nabla u)^{(p-1)}\big)+q(t)u^{(p-1)}=w(t)f(u),\\
\hskip 1.5cm u|_{\partial \Omega}=0.\end{array}\right.\
\end{eqnarray*}
This involves partial differential equations, but for symmetric
flows, the ordinary differential operator (perhaps in radial
form) is involved ( see, e.g., Binding \& Drabek \cite{bd}, del
Pino, Elgueta \& Manasevich \cite{dem}, del Pino \& Manasevich
\cite{dm}, Rabinowitz \cite{r} and Walter \cite{w}).
It is well-known that the problem $(BVP_1)$ has very similar
properties as the classical case when $p=2$, especially in the
one-dimensional case (cf. Erbe
\cite{er}, Kong \cite{k}, Lian, Wong \& Yeh \cite{lwy}, Naito \&
Tanaka \cite{nt} and the reference therein). It has been
investigated a good deal in the last twenty years or so under the
general heading of $p$-Laplacian.

In 2000, Erbe \cite{er} initiated the idea of connecting
the equation
\begin{equation}
\label{eq1.3}
-y''+q(t)y=w(t)f(y)
\end{equation}
under the separated boundary value conditions with the
Sturm-Liouville eigenvalue problem (SLEP). Using the fixed point
index method and comparing the values of $\frac{f(s)}{s},~s\in
(0,\infty)$, with the smallest eigenvalue of the corresponding
(SLEP), the existence of positive solutions of \eqref{eq1.3} was
established. But, due to the limitation of the approach, he only
discussed the case $q>0$ and certain boundary conditions and
nothing was found for the existence of solutions with zeros in
$(0,1)$.

Next, Naito \& Tanaka \cite{nt} compared the equation with the
$k$-th eigenvalue of the linear equation, and applied the method
of energy function and the Sturm-Picone comparison theorem to
establish the sufficient conditions for the existence of multiple
solutions with prescribed numbers of zeros in $(0,1)$ for the
case $q\equiv 0$, i.e., $$
\left\{\begin{array}{lll}
u''+w(t)f(u)=0~~~\hbox{on}~~~(0,1),\\
u(0)=u(1)=0.\\
\end{array}\right.$$

Recently, Kong \cite{k} generalized the results in \cite{nt} to
the case with nonzero $q$. This extension is not trivial due to
the fact that the energy function with a general $q$ may not be
nonnegative. He also obtained results on the nonexistence of
certain types of solutions. When $\frac{f(s)}{s}$,
$s\in(0,\infty)$, is between $f_0$ and $f_\infty$(see (C3) below
), the conditions for the existence and nonexistence become
necessary and sufficient.

In 2008, Naito \& Tanaka \cite{nt08} studied the Dirichlet
quasi-linear differential equation with $q\equiv 0$ by the
shooting method together with the qualitative theory. In the same
year, Lee \& Sim
\cite{ls08}  considered the same case, $q\equiv 0$, with
Dirichlet boundary conditions with a nonnegative measurable
function $w(x)$ on $(0,1)$ which may be singular at $x=0$ and/or
$x=1$. They gave the global analysis for sign-changing solutions
employing a bifurcation argument.

The aim of this paper is to establish the sufficient condition
for the existence of solutions of the BVP
\eqref{eq1.1}-\eqref{eq1.2} with prescribed numbers of zeros in
terms of the ratio $f(s)/s^{(p-1)}$ at infinity and zero.
Inspired by the ideas of  \cite{k,nt08,ls08}, we study  the BVP
\eqref{eq1.1}-\eqref{eq1.2} with a nonzero $q$ and  general separated boundary conditions.
 We  employ the generalized Pr\"{u}fer
substitution and comparison theorem in our arguments.
 Our results generalize partly those ones in
Kong \cite{k} and Naito \& Tanaka \cite{nt08}.

In this paper, we assume the following conditions hold:
\begin{enumerate}
\item[(C1)] $q$, $w\in C^1[0,1]$ and $w>0$ on $[0,1]$;
\item[(C2)] $f\in C(\mathbb{R})$, $f(s)>0$ for $s>0$,
$f(-s)=-f(s)$ for $s>0$, and $f$ is locally Lipschitz continuous
on $(0,\infty)$;
\item[(C3)] there exist $0\leq f_0, f_{\infty}\leq \infty$ such that $\lim_{s\rightarrow 0^+}\frac{f(s)}{s^{p-1}}=f_0$
and $\lim_{s\rightarrow\infty}\frac{f(s)}{s^{p-1}}=f_{\infty}$.
\end{enumerate}
The typical examples satisfying (C2)-(C3) are
\begin{enumerate}
\item[(i)] $f$ is a finite sum of sign power functions, $f(s)=\sum_{i=1}^nA_is^{(\ell_i)}$, with $A_i$, $\ell_i> 0$;
\item[(ii)] $f$ is an odd exponential function such as $f(s)=A\sinh Bs$, where $A$, $B> 0$.
\end{enumerate}

In order to discuss our results, we  compare the BVP \eqref{eq1.1} -
\eqref{eq1.2} with the $p$-Laplacian eigenvalue problem
\begin{equation}
\label{eq1.4}
-(y'^{(p-1)})'=(p-1)(\lambda w(x)-q(x))y^{(p-1)} \mbox{ on }
(0,1),
\end{equation}
coupled with the boundary condition \eqref{eq1.2}. It is known
that BVP
\eqref{eq1.4} and \eqref{eq1.2} has a countable number of eigenvalues $\lambda_i$ satisfying
$$-\infty<\lambda_0<\lambda_1<\lambda_2<\cdots
\cdots< \lambda_n< \cdots\ \to\infty\ ,$$
and the eigenfunction corresponding to $\la_{n}$ has $n$ zeros in
$(0,1)$ (cf. Binding
\& Drabek
\cite{bd}, Reichel \& Walter \cite{rw99}). The followings
are our main results.

\newtheorem{th1.1}{Theorem}[section]
\begin{th1.1}
\label{th1.1}
\begin{enumerate}
\item[(i)] For all $y\in (0,\infty)$, if $\frac{f(y)}{y^{p-1}}< \lambda_n$ for some $n$,
then BVP \eqref{eq1.1} -\eqref{eq1.2} has no solution with exactly $i$ zeros in $(0,1)$
for any $i\geq n$.
\item[(ii)] For all $y\in(0,\infty)$, if $\lambda_n<\frac{f(y)}{y^{p-1}}$ for some $n$,
then BVP \eqref{eq1.1} -\eqref{eq1.2} has no solution with exactly $i$ zeros in $(0,1)$ for any $i\leq n$.
\item[(iii)] For all $y\in (0,\infty)$, if $\frac{f(y)}{y^{p-1}}\neq \lambda_n$ for any $n$,
then BVP \eqref{eq1.1} -\eqref{eq1.2} has no nontrivial solution.
\end{enumerate}
\end{th1.1}

\newtheorem{th1.2}[th1.1]{Theorem}
\begin{th1.2}
\label{th1.2} Assume there exists $n$ such that
either $\lambda_n\in (f_0, f_{\infty})$ or $(f_{\infty}, f_0)$.
Then BVP \eqref{eq1.1} -\eqref{eq1.2} has a solution with exactly
n zeros in $(0,1)$.
\end{th1.2}

The combination of Theorems \ref{th1.1} and \ref{th1.2} leads to
the following.

\newtheorem{th1.3}[th1.1]{Corollary}
\begin{th1.3}
\label{th1.3}
  Assume $\frac{f(y)}{y^{p-1}}\in (f_0,f_{\infty})$ or $(f_{\infty},f_0)$ for all $y\in (0,\infty)$. Then for some $n$ such that $\la_{n}>0$,
BVP \eqref{eq1.1} -\eqref{eq1.2} has a solution with exactly n
zeros in $(0,1)$ if and only if $\la_{n}\in (f_0,f_{\infty})$ or
$(f_{\infty},f_0)$.

\end{th1.3}

This paper is organized as follows. After the introduction, we
establish the global existence and uniqueness of the solutions of
the initial value problems associated with \eqref{eq2.1} in
Section 2. In Section 3, we give some technical lemmas for the
proof of Theorem \ref{th1.2}. We give the proofs of the
main theorems in Section 4. In the appendix, we give the proof of Proposition \ref{th2.2} represented in Section 2.

\section{Results on initial value problems}
\setcounter{equation}{0} \hskip0.25in We first establish the basic properties of the system
consisting of \eqref{eq1.1} with the initial condition
\begin{equation}
\label{eq2.1}
y(0)=\eta_{1},~~~y'(0)=\eta_{2},
\end{equation}
where $\eta_{1},\eta_{2}\in \mathbb{R}$ is a parameter. It is
known the solution of this IVP, if it exists, may not be unique
under a general condition (cf.\ Walter
\cite[p.181]{w} and Reichel \& Walter \cite[Theorem
4]{rw} for the radial case). Now we show that the global existence
and uniqueness are guaranteed under the conditions (C1)-(C3).
First motivated by \cite{k}, we introduce a generlized energy
function $E[y](x)$ to derive the global existence.

\newtheorem{th2.1}{Proposition}[section]
\begin{th2.1}
\label{th2.1}
For any $\eta_{1},\eta_{2}\in \mathbb{R}$, the IVP \eqref{eq1.1}
and \eqref{eq2.1} has a solution $y(x)$ which exists over the whole
interval $[0,1]$.
\end{th2.1}
\begin{proof}
Note that the IVP \eqref{eq1.1} and \eqref{eq2.1} can be written
by
\begin{equation}
\label{eq2.2}
\left\{\begin{array}{l}
y'(x)=z(x)^{(p^*-1)},\\
z'(x)=(p-1)q(x)y(x)^{(p-1)}-(p-1)w(x)f(y(x)),
\end{array}\right.\
\end{equation}
with $y(0)=\eta_{1}$ and $z(0)=\eta_{2}^{(p-1)}$ where
$p^*=\frac{p}{p-1}$ is the conjugate exponent of $p$. Then the
local existence of a solution is guaranteed by  Peano existence
theorem. We will divide the proof into two steps to prove the
global existence.

\begin{enumerate}
\item[(i)] Consider the case $f_{\infty}=\infty$. Suppose $y(x)$ does not exist on the whole interval $[0,1]$.
Without loss of generality, we may assume $y(x)$ exists on a
maximal right interval $[0,c)$ for some $c\in (0,1)$. Then $y(x)$
is unbounded on $[0,c)$. For otherwise, if $y(x)$ is bounded on
$[0,c)$, then, from
\eqref{eq1.1}, $\lim_{x\to c }y'(x)$ is bounded. This implies
$y(x)$ can be extended through $c$. This contradicts that $[0,c)$
is the maximal right interval of existence for $y(x)$. Now define
the generalized energy function for the solution $y$ as
    \begin{equation}
    \label{eq2.3}
    E[y](x)=\frac{|y'(x)|^p}{p}-\frac{1}{p}q(x)|y(x)|^p+w(x)F(y(x)),
    \end{equation}
where $F(y)=\int_0^yf(s)ds$. Then \eqref{eq1.1} implies that
$$E[y]'(x)=-\frac{1}{p}q'(x)|y(x)|^p+\frac{w'(x)}{w(x)}w(x)F(y(x)).$$
Let $k=\max \{\frac{|w'(x)|}{w(x)}:x\in [0,1]\}$. Then
    \begin{equation}
    \label{eq2.4}
    E[y]'(x)\leq -\frac{k+1}{p}q(x)|y(x)|^p+\frac{1}{p}[(k+1)q(x)-q'(x)]|y(x)|^p+kw(x)F(y(x)).
    \end{equation}
    Because $w> 0$ is continuous and $q$, $q'$ are bounded on $[0,1]$, we can find a constant $h> 0$
    such that
    \begin{equation}
    \label{eq2.5}
     \frac{h}{p}[(k+1)q(x)-q'(x)]\leq w(x)~~and~~\frac{h}{p}|q(x)|\leq w(x).
     \end{equation}
 Since $f_{\infty}=\infty$, there exists $M> 0$ such that $|y|^p\leq hF(y)$ for $|y|\geq M$.
    Define $I_1=\{x\in [0,c):|y(x)|\leq M\}$ and $I_2=\{x\in [0,c):|y(x)|> M\}$.
    Then from \eqref{eq2.4}, there exists $N> 0$ such that $E'[y](x)\leq N$ for $x\in I_1$, and
    for $x\in I_2$ $$E[y]'(x)\leq (k+1)[-\frac{1}{p}q(x)|y(x)|^p+w(x)F(y(x))]\leq (k+1)E[y](x).$$
    Hence, from the second inequality in \eqref{eq2.5} that $E[y](x)\geq 0$ on $I_2$; we have $$E[y]'(x)\leq N+(k+1)E[y](x),~~~x\in [0,c).$$
    Integrating both sides of the above inequality, we get that for $x\in [0,c)$,
    $$E[y](x)\leq E[y](0)+N+\int_0^x(k+1)E[y](t)dt.$$ By the Gronwall inequality,
    $$E[y](x)\leq (E[y](0)+N)\exp[(k+1)x],~~~for~~x\in [0,c).$$ Therefore,
    \begin{equation}
    \label{eq2.6}
    \limsup_{x\rightarrow c-}E[y](x)< \infty.
    \end{equation}
    On the other hand, since $y(x)$ is unbounded on $[0,c)$, there exists a sequence $t_n\rightarrow c-$
    such that $|y(t_n)|\rightarrow \infty$. Hence $\lim_{n\rightarrow \infty}\frac{F(y(t_n))}{|y(t_n)|^p}=\frac{f_{\infty}}{p}=\infty$.
    By \eqref{eq2.3}, $$E[y](t_n)\geq (-\frac{1}{p}q(t_n)+w(t_n)\frac{F(y(t_n))}{|y(t_n)|^p})|y(t_n)|^p\rightarrow \infty~~as~~n\rightarrow \infty.$$
    This contradicts with \eqref{eq2.6}.
\item[(ii)] Let $f_{\infty}< \infty$. Integrating \eqref{eq2.2} over any compact interval $[0,k]\subset [0,\infty)$, we obtain
\begin{eqnarray}
y(x)&=&\eta_1+\int_0^xz(t)^{(p^*-1)}dt,\label{eq2.7}\\
z(x)&=&\eta_2^{(p-1)}+(p-1)\int_0^xq(t)y(t)^{(p-1)}dt-(p-1)\int_0^xw(t)f(y(t))dt,\label{eq2.8}
\end{eqnarray}
for any $x\in [0,k]$. Note that $p^*-1=\frac{1}{p-1}$. Since
$f_{\infty}< \infty$, we have $$|f(y)|\leq
c_1|y|^{p-1}~~~\hbox{for}~~~|y|\geq M',$$ where $c_1$ and $M'$ are
some positive numbers. For $|y|<M'$, it is easy to obtain the
boundedness of $z(x)$ by \eqref{eq2.8}. So, for $|y|\geq M'$ and
$x\in [0,k]$, it follows from H\"{o}lder inequality that
$$|z(x)|\leq c_2+c_3(\int_0^x|y(t)|^pdt)^{\frac{p-1}{p}}.$$ i.e.,
$$|z(x)|^{\frac{p}{p-1}}\leq
c_2^{\frac{p}{p-1}}+(\int_0^x|y(t)|^pdt)(c_3^{\frac{p}{p-1}}+O(1)).$$
Thus
\begin{equation}
\label{eq2.9} |z(x)|^{\frac{p}{p-1}}\leq
c_4+c_5\int_0^x|y(t)|^pdt,~~~for~~x\in [0,k].
\end{equation}
Similar arguments, we can obtain
\begin{equation}
\label{eq2.10}
|y(x)|^p\leq c_6+c_7\int_0^x|z(t)|^{p^*}dt,,~~~for~~x\in [0,k].
\end{equation}

 From \eqref{eq2.9} and
\eqref{eq2.10}, we have $$|y(x)|^p+|z(x)|^{p^*}\leq
c(k;p)+d(k;p)\int_0^x(|y(t)|^p+|z(t)|^{p^*})dt,$$ where $c(k;p)$
and $d(k;p)$ are some positive constants depending on $p$ and $k$.
By Gronwall inequality, $$|y(x)|^p+|z(x)|^{p^*}\leq
c(k;p)\exp[d(k;p)x]<
\infty.$$
\end{enumerate}
Therefore, the solution exists over the whole interval $[0,1]$.
\end{proof}

Following the ideas of \cite{w}, \cite{k} and \cite{nt08}, we can
prove the following uniqueness of the solution of the IVP, which
will be proven in the appendix.

\newtheorem{th2.2}[th2.1]{Proposition}
\begin{th2.2}
\label{th2.2} For any $\eta_{1},\eta_{2} \in\mathbb{R}$, the solution $y(x ; \eta_{1}, \eta_{2})$ of the IVP
\eqref{eq1.1} and \eqref{eq2.1} is unique on $[0,1]$. Furthermore, $y(x ; \eta_{1},
\eta_{2})$ and $y'(x ; \eta_{1}, \eta_{2})$ are continuous in $(x;\eta_{1}, \eta_{2} )\in [0,1]\times
\mathbb{R}^{2}$.
\end{th2.2}

\section{Some technical lemmas}
\setcounter{equation}{0}

\hskip0.25in In this section, we will derive three lemmas for the proof of main theorems. First we consider the IVP consisting \eqref{eq1.1} with the initial condition
\begin{equation}
\label{eq3.1} y(0)=\rho S_p(\alpha),~~~y'(0)=\rho S_p'(\alpha),
\end{equation}
where $\rho> 0$ is a parameter. Denote by $y(x;\rho)$ the solution
of
\eqref{eq1.1} and \eqref{eq3.1}. Consider the modifier
Pr\"{u}fer substitution $$y(x;\rho) =
r(x;\rho)S_p(\theta(x;\rho))\ ,\qquad y'(x;\rho) =
r(x;\rho)S'_p(\theta(x;\rho))\ .$$ Then we have
$\theta(0;\rho)=\alpha, r(x;\rho)=(|y(x;\rho)|^p+|y'(x;\rho)|^p)^{1/p}> 0$
and
\begin{eqnarray}
\theta'(x;\rho)&=&|S_p'(\theta(x;\rho))|^p+\frac{w(x)f(y(x;\rho))S_p(\theta(x;\rho))}{r(x;\rho)^{p-1}}-q(x)|S_p(\theta(x;\rho)|^p\label{eq3.2}\
,\\
r'(x;\rho)&=&S_p'(\theta(x;\rho))\left[(1+q(x))r(x;\rho)S_p^{(p-1)}(\theta(x;\rho))-\frac{w(x)f(y(x;\rho))}{r(x;\rho)^{p-2}}\right]\label{eq3.3}\
.
\end{eqnarray}
Similarly, the Pr\"{u}fer angle $\phi_n$ for
\eqref{eq1.4} and \eqref{eq3.1} with $\la=\la_{n}$ satisfies
\begin{equation}
\label{eq3.4}
\phi_n'(x;\rho)=|S_p'(\phi_n(x;\rho))|^p+[\lambda_nw(x)-q(x)]|S_p(\phi_n(x;\rho))|^p\
.
\end{equation}

\newtheorem{th3.1}{Lemma}[section]
\begin{th3.1}
\label{th3.1}
\begin{enumerate}
\item[(a)]
Assume $f_0< \lambda_n$ for some $n$. Then there exists a
sufficiently small $\rho_* $ such that $\theta(1;\rho)<
n\pi_p+\beta$ for all $\rho
\in(0,\rho_*)$.
\item[(b)]
Assume $f_0> \lambda_n$ for some $n$. Then there exists  a
sufficiently small $\rho_* $ such that $\theta(1;\rho)>
n\pi_p+\beta$ for all $\rho
\in(0,\rho_*)$.
\end{enumerate}

\end{th3.1}
\begin{proof}

We  give the proof of (a) here. The proof of  part (b) is similar
and will be omitted.

Since $f_0< \lambda_n< \infty$, $\frac{f(y)}{y^{(p-1)}}$ can be
continuously extended to $y=0$ and there exists $\delta >0$ such
that $$\frac{f(y)}{y^{(p-1)}}< \lambda_n~~~for~~~|y|< \delta.$$
Since $y\equiv 0$ is a solution of \eqref{eq1.1}, by the
continuous dependence of solutions on the initial conditions,
there exists $\rho_*> 0$ such that $|y(x;\rho)|< \delta$ for
$\rho < \rho_*$ and $x\in [0,1]$. From \eqref{eq3.2}, for
$\rho<\rho_*$ and $x\in[0,1]$, we have
$$\theta'(x;\rho)<|S_p'(\theta(x;\rho))|^p+[\lambda_nw(x)-q(x)]|S_p(\theta(x;\rho))|^p.$$
Let $u_n(x;\rho)$ be the solution of the IVP
\eqref{eq1.4} and \eqref{eq3.1} with $\lambda=\lambda_n$ and let $\phi_n$ be its Pr\"{u}fer
angle. Then $u_n(x;\rho)$ is an eigenfunction of the BVP
\eqref{eq1.4} and \eqref{eq1.2} corresponding to the eigenvalue
$\lambda=\lambda_n$; thus $\phi_n(1;\rho)=n\pi_p+\beta$. By the
comparison theorem, we obtain that $\theta(1;\rho)<
\phi_n(1;\rho)$. This completes the proof.
\end{proof}

\newtheorem{th3.3}[th3.1]{Lemma}
\begin{th3.3}
\label{th3.3}
For $M> 0$ and $\rho > 0$ define $I_{M,\rho}=\{x\in
[0,1]:|y(x;\rho)|< M\}$. Then for any $M$, $L> 0$, there exists a
sufficiently large $\rho^* >0$ such that $|y'(x;\rho)|> L$ for
$\rho
>\rho^*$ and $x\in I_{M,\rho}$.
\end{th3.3}
\begin{proof}
\begin{enumerate}
\item[(i)] Let $f_{\infty}< \infty$. For $\rho >0$ and from \eqref{eq3.3}, it is easy to find $K_1>0$ such that
    \begin{equation}
    \label{eq3.8}
    r'(x;\rho)\geq -K_1~~~for~~x\in I_{M,\rho}.
    \end{equation}
    Since $f_{\infty}< \infty$, there exists $K_2>0$ such that $|\frac{f(y(x;\rho))}{y(x;\rho)^{(p-1)}}|\leq K_2$ for $x \in[0,1]\backslash I_{M,\rho}$.
By \eqref{eq3.3}, we have, for $x \in[0,1]\backslash I_{M,\rho}$,
    \begin{eqnarray}
    r'(x;\rho)&=&r(x;\rho)S_p'(\theta(x;\rho))S_p^{(p-1)}(\theta(x;\rho))
    \left[1+q(x)-\frac{ w(x)f(y(x;\rho))}{y^{(p-1)}(x;\rho)}\right]\nonumber\\
&\geq&  -r(x;\rho)[|1+q(x)|+K_2w(x)]\nonumber\\
&\geq& -K_3r(x;\rho)\ ,
        \label{eq3.9}
     \end{eqnarray}
where $K_3= \max \{|1+q(x)|+K_2w(x) : x\in [0,1]\}$.
    Combining \eqref{eq3.8} and \eqref{eq3.9}, $$r'(x;\rho)\geq -K_1-K_3r(x;\rho)~~~for~~\rho>0~~and~~x\in[0,1].$$
      Solving the above linear differential inequality for $x\in [0,1]$,  we have $$r(x;\rho)\geq -\frac{K_1}{K_3}+(\rho+\frac{K_1}{K_3})e^{-K_3x}\rightarrow \infty~~~as~~\rho\rightarrow \infty$$
     uniformly in $[0,1]$.
    Therefore, for any $L> 0$, there exists $\rho^*> 0$ such that $\rho >\rho^*$
    and $x\in I_{M,\rho}$ $$M^p+|y'(x;\rho)|^p \geq r(x;\rho)> M^p+L^p.$$ This leads to have that $|y'(x;\rho)|> L$.

\item[(ii)] Let $f_{\infty}=\infty$. Recall the generalized energy
function $E[y](x;\rho)$ defined as \eqref{eq2.3},
\begin{equation}
\label{eq3.45}
E[y](x;\rho)=\frac{|y'(x;\rho)|^p}{p}-\frac{1}{p}q(x)|y(x;\rho)|^p+w(x)F(y(x;\rho))
\end{equation}
where $F(y)=\int_{0}^{y}f(s)ds$.
 Then, letting $k=\max\{\frac{|w'(x)|}{w(x)}: x\in[0,1]\}$,
    \begin{eqnarray}
    E[y]'(x;\rho)&=&-\frac{1}{p}q'(x)|y(x;\rho)|^p+w'(x)F(y(x;\rho))\nonumber\\
    &\geq&\frac{k+1}{p}q(x)|y(x;\rho)|^p-\frac{1}{p}[(k+1)q(x)+q'(x)]|y(x;\rho)|^p\nonumber\\
    &~~~&-kw(x)F(y(x;\rho)).\label{eq3.5}
    \end{eqnarray}
Because $w> 0$ is continuous and $q$, $q'$ are bounded on
$[0,1]$, we can find a constant $h> 0$
    such that
    \begin{equation}
    \label{eq4.10}
     \frac{h}{p}\left|(k+1)q(x)+q'(x)\right| \leq w(x)~~\mbox{and}~~\frac{h}{p}|q(x)|\leq w(x).
     \end{equation}
    Since $f_{\infty}=\infty$, when $M$ is sufficiently large, we have $|y|^p\leq hF(y)$ for all $|y|\geq M$.
    Then from \eqref{eq3.5}, when $x\in I_{M,\rho}$, $E[y]'(x;\rho)$ is bounded below.
    But  for $x\in [0,1]\backslash I_{M,\rho}$, by (\ref{eq4.10}), \
    $E[y]'(x;\rho)\geq -(k+1)E[y](x;\rho)$.
      Hence, for all $x\in[0,1]$,
      $$E[y]'(x;\rho)\geq -N-(k+1)E[y](x;\rho).$$
      Solving the above linear differential inequality for $x\in [0,1]$, we obtain
    \begin{equation}
    \label{eq3.6}
    E[y](x;\rho)\geq -\frac{N}{k+1}+(E[y](0;\rho)+\frac{N}{k+1})\exp[-(k+1)x].
    \end{equation}
    Note that, from the initial condition \eqref{eq3.1},  $$E[y](0;\rho)=\rho^p[\frac{|S_p'(\alpha)|^p}{p}-\frac{q(x)|S_p(\alpha)|^p}{p}+w(x)\frac{F(\rho S_p(\alpha))}{\rho^p}].$$
     When $\alpha=0$, $$\lim_{\rho\rightarrow \infty}E[y](0;\rho)=\lim_{\rho\rightarrow \infty}\frac{\rho^p}{p}=\infty,$$
    and when $\alpha \in(0,\pi_p)$, it follows from $f_{\infty}=\infty$
    that
    $$\lim_{\rho\rightarrow \infty}\frac{F(\rho S_p(\alpha))}{\rho^p}=\lim_{y\rightarrow \infty}\frac{F(y)}{y^p}|S_p(\alpha)|^p=\infty.$$
     So we have $\lim_{\rho\rightarrow \infty}E[y](0;\rho)=\infty$. Therefore, from \eqref{eq3.6} we get that
    \begin{equation}
    \label{eq3.7}
    \lim_{\rho\rightarrow \infty}E[y](x;\rho)=\infty~~~uniformly~~for~~x\in[0,1].
    \end{equation}
    Note that, in \eqref{eq3.45}, the term  $|-\frac{1}{p}q(x)|y(x;\rho)|^p+w(x)F(y(x;\rho))|$ is uniformly bounded for all $\rho > 0$ and $x\in I_{M,\rho}$.
    So, by \eqref{eq3.45} and \eqref{eq3.7}, we may choose $\rho^*$ such that  $\rho>\rho^*$ and $x\in I_{M,\rho}$, $|y'(x;\rho)|> L$.

\end{enumerate}
\end{proof}

\newtheorem{th3.4}[th3.1]{Lemma}
\begin{th3.4}
\label{th3.4}
\begin{enumerate}
  \item[(a)]Assume $f_{\infty}>\lambda_n$ for some $n$. Then there
  exists a
sufficiently large $\rho^*$ such that $\theta(1;\rho)>
n\pi_p+\beta$ for all $\rho
\in(\rho^*,\infty)$.
  \item[(b)]Assume $f_{\infty}< \lambda_n$ for some $n\in \mathbb{N}_k$. Then there exists a
sufficiently large $\rho^*$
 such that $\theta(1;\rho)< n\pi_p+\beta$ for all $\rho \in(\rho^*,\infty)$.

\end{enumerate}

\end{th3.4}
\begin{proof}

We  give the proof of (a) here. The proof of  part (b) is similar.

Assume the contrary. Then there exists ${\rho_l}$ with
$\rho_l\rightarrow \infty$ such that $\theta(1;\rho_l)\leq
n\pi_p+\beta$. This implies that $y(x;\rho_l)$ has at most $n$ zeros
in $(0,1)$. Since $f_{\infty}> \lambda_n$, we can choose $\lambda > 0$ such that $\lambda_n<
\lambda < f_{\infty}$ and take $M> 0$ so that
$$\frac{f(y(x;\rho))}{y(x;\rho)^{(p-1)}}\geq
\lambda~~~for~~|y(x;\rho)|\geq M.$$ We divide the proof into several
steps:
\begin{enumerate}
\item[(i)] We claim that the measure of $I_{M,\rho}$ tends to zero as $\rho=\rho_l\rightarrow \infty$.
It is easy to see that for each $\rho=\rho_l$, $I_{M,\rho}\cap
(0,1)$ is an open set and hence is a union of disjoint intervals
in $(0,1)$, i.e., $$I_{M,\rho}\cap (0,1)=\cup_{i=1}^j(a_i,b_i),$$
where $0\leq a_i< b_i\leq 1$. If $0< a_i$ and $b_i<1$, by Lemma
\ref{th3.3}, for $\rho=\rho_l$ sufficiently large, $y(x;\rho)$ is
monotone on $[a_i,b_i]$, and hence
$|y(a_i;\rho)|=|y(b_i;\rho)|=M$ and $y(a_i;\rho)y(b_i;\rho)< 0$.
This implies that $(a_i,b_i)$ contains exactly one zero of
$y(x;\rho)$, so $j\leq n+2$. Applying Lemma \ref{th3.3} again, for
any $L>0$ there exists $\rho(L)>0$ such that if
$\rho=\rho_l>\rho(L)$, then $y'(x;\rho)$ has the same sign and
$|y'(x;\rho)|>L$ in $(a_i,b_i)$ for each $1\leq i\leq j$. Thus,
$$2M= |y(b_i;\rho)-y(a_i;\rho)|=|\int_{a_i}^{b_i}y'(t;\rho)dt|>
L(b_i-a_i).$$ i.e., $b_i-a_i< \frac{2M}{L}$. So,
$\|I_{M,\rho}\|\leq \frac{2(n+2)M}{L}$, where $\|\cdot\|$ is the
Lebesgue measure. Therefore,
    \begin{equation}
    \label{eq3.10}
    \lim_{\rho_l\rightarrow \infty}\|I_{M,\rho_l}\|=0.
    \end{equation}

\item[(ii)] Next, we try to reach a contradiction to our assumption. For each $\rho=\rho_l$,
 let $\phi(x;\rho)$  and $\phi_n(x;\rho)$ be the Pr\"{u}fer angles
of the solution of \eqref{eq1.4} and \eqref{eq3.1} with $\lambda$ and
$\lambda_n$, respectively. Then $\phi_n(1;\rho)=n\pi_p+\beta$ and
hence, by the comparison theorem,
$\phi(1;\rho)=n\pi_p+\beta+\epsilon$ for some $\epsilon> 0$.
Recall that $\phi(x;\rho)$ satisfies
    \begin{equation}
    \label{eq3.11}
    \phi'(x;\rho)=|S_p'(\phi_(x;\rho))|^p+[\lambda w(x)-q(x)]|S_p(\phi(x;\rho))|^p\equiv G(x,\rho,\phi),
    \end{equation}
    and $\phi(0;\rho)=\alpha$.
On the other hand, define
    \begin{eqnarray*}g(x;\rho)=\left\{\begin{array}{l}\frac{f(y(x;\rho))}{y(x;\rho)^{(p-1)}},~~|y(x;\rho)|<M\\
    \lambda,~~~~~~~~~~|y(x;\rho)|\geq M.
    \end{array}\right.\
    \end{eqnarray*}
and let $\theta(x;\rho)$ be the Pr\"{u}fer angle of
\eqref{eq1.1} and \eqref{eq3.1}. Then,
    by \eqref{eq3.2},$$\theta'(x;\rho)\geq |S_p'(\theta(x;\rho))|^p+w(x)g(x;\rho)|S_p(\theta(x;\rho))|^p-q(x)|S_p(\theta(x;\rho))|^p\equiv F(x,\rho,\theta).$$
     Let $\bar{\theta}(x;\rho)$ be the solution of the equation,
    \begin{equation}
    \label{eq3.12}
    \bar{\theta}'(x;\rho)=F(x,\rho,\bar{\theta})
    \end{equation}
    satisfying $\bar{\theta}(0;\rho)=\alpha$. From \eqref{eq3.11} and \eqref{eq3.12} we have that for $\rho=\rho_l$ and $x\in [0,1]$, \begin{eqnarray}
    \bar{\theta}(x;\rho)-\phi(x;\rho)&=&\int_0^x(F(s,\rho,\bar{\theta})-G(s,\rho,\phi))ds\nonumber\\
    &=&\int_0^x[(F(s,\rho,\bar{\theta})-G(s,\rho,\bar{\theta}))+(G(s,\rho,\bar{\theta})-G(s,\rho,\phi))]ds\nonumber\\
    &=&\int_0^xw(s)[g(s;\rho)-\lambda]|S_p(\bar{\theta}(s;\rho))|^pds\nonumber\\
    &~~~&+\int_0^x\frac{\partial}{\partial \phi}G(s,\rho,\zeta)[\bar{\theta}(s;\rho)-\phi(s;\rho)]ds\label{eq3.13}
    \end{eqnarray}
    where $\zeta(s;\rho)$ is between $\bar{\theta}(s;\rho)$ and $\phi(s;\rho)$. Since $g(x;\rho)=\lambda$ for $x\in [0,1]\backslash I_{M,\rho}$ and
    $g(x;\rho)$ is continuous on $I_{M,\rho}$,  we have, by
    \eqref{eq3.10},
    \begin{equation}
    \label{eq3.14}
    |\int_0^xw(s)[g(s;\rho)-\lambda]|S_p(\bar{\theta}(s;\rho))|^pds|\leq \int_{I_{M,\rho}}w(s)|g(s;\rho)-\lambda|ds \rightarrow 0
    \end{equation}
    as $\rho=\rho_l \rightarrow \infty$. Note that $|\frac{\partial}{\partial \phi}G(x,\rho,\phi)|$ is uniformly bounded by some $K> 0$
    for all $x\in [0,1]$. Thus, by \eqref{eq3.13} and \eqref{eq3.14}, for any $\delta > 0$ there exists a large $\rho^*$ such that for $\rho\in (\rho^*,\infty)$,
    $$|\bar{\theta}(x;\rho)-\phi(x;\rho)|< \delta +\int_0^xK|\bar{\theta}(s;\rho)-\phi(s;\rho)|ds.$$
    By Gronwall inequality, we have $$|\bar{\theta}(x;\rho)-\phi(x;\rho)|< \delta e^{Kx} < \epsilon$$
    if $\delta <\epsilon e^{-K}$. Hence, $\bar{\theta}(x;\rho)> \phi(x;\rho)-\epsilon$
    on $[0,1]$. Furthermore, $$\theta(1;\rho)\geq \bar{\theta}(1;\rho)> \phi(1;\rho)-\epsilon=n\pi_p+\beta.$$ This contradicts with our assumption. So the proof is completed.
\end{enumerate}
\end{proof}

\section{Proofs of Theorem \ref{th1.1} \& \ref{th1.2}}
\setcounter{equation}{0}

\begin{proof}[Proof of Theorem \ref{th1.1}]
\begin{enumerate}
\item[(i)] Assume the contrary that there exists a solution $y(x)$ of \eqref{eq1.1}-\eqref{eq1.2} with exactly $i$ zeros in $(0,1)$
 for some $i\geq n$. Let $\bar{w}(x)=w(x)\frac{f(y(x))}{y(x)^{(p-1)}}$. Then $\bar{w}(x)$ is continuous on $[0,1]$ since $f_0< \infty$.
  Denote that $\theta(x)$ is the Pr\"{u}fer angle of $y(x)$. Then $\theta(x)$ satisfies \eqref{eq3.2} and $\theta(1)=(i+1)\pi_p$, $i\geq n$.
   Note that $\frac{f(y(x))}{y(x)^{(p-1)}}\neq \infty$ in this case, and
    $$w(x)\frac{f(y(x))S_p(\theta(x))}{r(x)^{(p-1)}}=\frac{w(x)f(y(x))}{y(x)^{(p-1)}}|S_p(\theta(x))|^p< \lambda_nw(x)|S_p(\theta(x))|^p.$$
     By the comparison theorem, we find that $\theta(1)<\phi_n(1)=(n+1)\pi_p$. This reaches a contradiction.

\item[(ii)] Assume the contrary that \eqref{eq1.1}-\eqref{eq1.2} has a solution $y(x)$ with exactly $i$ zeros in $(0,1)$ for some $i\leq n$.
By similar argument as the above, we have that $\theta'(x)>
\phi'_n(x)$ a.e. on $[0,1]$. By the comparison theory,  we obtain that $\theta(1)> \phi_n(1)=(n+1)\pi_p$.
\item[(iii)] The assumption implies that either
\begin{enumerate}
\item[(a)] $\lambda_n< \frac{f(y(x))}{y(x)^{p-1}}< \lambda_{n+1}$ for some $n$ for all $y\in (0,\infty)$

or

\item[(b)] $0< \frac{f(y(x))}{y(x)^{p-1}}<\lambda_0$ for all $y\in (0,\infty)$ if $k=0$.
\end{enumerate}
By the uniqueness, we have that the number of zeros of $y(x)$ is
finite. Then the conclusion follows from $(i)$ and $(ii)$
immediately.
\end{enumerate}
\end{proof}

\begin{proof}[Proof of Theorem \ref{th1.2}] Here we consider the case of $f_0< \lambda_n<f_{\infty}$.
By Lemma \ref{th3.1}, there exists $\rho_*>0$ such that
$\theta(1;\rho)<n\pi_p+\beta$ for all $\rho\in (0,\rho_*)$. By
Lemma \ref{th3.4}, there exists $\rho^*> \rho_*$ such that
$\theta(1;\rho)> n\pi_p+\beta$ for all $\rho
\in(\rho^*,\infty)$. By continuous dependence on parameters, there exists $\rho_n\in [\rho_*,\rho^*]$ such that
$\theta(1;\rho_n)=n\pi_p+\beta$. This implies that $y(x;\rho_n)$
is a solution of the BVP \eqref{eq1.1}-\eqref{eq1.2} with exactly
$n$ zeros in $(0,1)$.

The other case is similar by applying Lemma \ref{th3.1}(b) and
\ref{th3.4} (b).

\end{proof}
\section*{Appendix. The proof of Proposition \ref{th2.2}}
In the proof, we need the following lemma (\cite{w}, p180).
\newtheorem*{th5.1}{Lemma A}
\begin{th5.1}
Let $W\in C^1(I)$, $x_0\in I$ and $W(x_0)=0$, where $I$ is a compact interval containing $x_0$. Denote $\|W\|_x$ the maximum of $W$ in the interval from $x_0$ to $x$. Then
\begin{equation}
\label{eq5.1}
|W'(x)|\leq K\| W\|_x~~~in~I~implies~~W=0~~for~~|x-x_0|\leq \frac{1}{K},~x\in I.
\end{equation}
\end{th5.1}
Now, it suffices to show the uniqueness of a local solution of the IVP. We will divide the proof into the following cases.
\begin{enumerate}
\item[(i)] Let $\eta_1\cdot \eta_2\neq 0$. Since the right hand side of \eqref{eq2.2} is locally Lipschitz continuous in $y$, $z\in \mathbb{R}\backslash\{0\}$, the local solution of \eqref{eq1.1} and \eqref{eq2.1} is unique for this case.
\item[(ii)] Let $\eta_1=\eta_2=0$. In this case we apply two versions of energy functions to prove the solution $y(x;\eta_1,\eta_2)\equiv 0$.
    \begin{enumerate}
    \item[(a)] Let $f_0=\infty$. Let $E[y](x)$ be defined as \eqref{eq2.3}. Similar arguments in the proof of Proposition \ref{th2.1},
     we have
        \begin{equation}
        \label{eq5.2}
        E[y]'(x)\leq -\frac{(k+1)}{p}q(x)|y(x)|^p+\frac{1}{p}[(k+1)q(x)-q'(x)]|y(x)|^p+kw(x)F(y(x)),
        \end{equation}
       and we can choose $h> 0$ such that
        \begin{equation}
        \label{eq5.3}
        \frac{h}{p}|(k+1)q(x)-q'(x)|\leq w(x)~and~\frac{h}{p}|q(x)|\leq w(x)~~on~~[0,1].
        \end{equation}
        Since $f_0=\infty$, there exists $\delta >0$ such that
        \begin{equation}
        \label{eq5.4}
        |y(x)|^p< hF(y(x))~~ for~~ |y(x)|<\delta.
         \end{equation}
         In this case there is a subinterval $[0,c)$, where $c\in (0,1]$, such that $|y(x)|<\delta$ on $[0,c)$.
        From \eqref{eq5.2} and \eqref{eq5.3}, for $x\in[0,c)$ we have
        $$E[y]'(x)\leq (k+1)E[y](x).$$i.e., $$E[y](x)\leq E[y](0)e^{(k+1)x}~~on~~[0,c).$$
        But $E[y](0)=0$, and so $E[y](x)\leq 0$ for $x\in[0,c)$. In fact, if there exists $x_1\in (0,c)$ such that $y(x_1)\neq 0$, then by \eqref{eq5.3} and \eqref{eq5.4}
        we have $E[y](x_1)> 0$. This leads a contradiction.
       \item[(b)] Let $f_0< \infty$. For the above $\delta>0$, there exists $d_1>0$ such that $|f(y)|\leq d_1|y(x)|^{p-1}$ for $|y(x)|\leq \delta$.
         Note that $|y(x)|< \delta$ on $[0,c)$.
         Applying the similar arguments in the proof of Proposition \ref{th2.1}, by \eqref{eq2.7}-\eqref{eq2.8}, we have
         $$|y(x)|^p+|z(x)|^{p^*}\leq
         d(c) \int_0^x(|y(t)|^p+|z(t)|^{p^*})dt,$$where $d(c)$ is
         some positive constant. By Gronwall inequality, it completes this
         case.
    \end{enumerate}
\item[(iii)] Let $\eta_1=0$ and $\eta_2\neq 0$. Without loss of generality, say $\eta_2=1$. Then
\begin{equation}
\label{eq5.5}
\frac{1}{2}|x-0|< |y(x)|<2|x-0|~~near~~0.
\end{equation}
Now we assume that $y_1(x)$ and $y_2(x)$ are two local solution of the IVP with the same initial condition.
Then, by \eqref{eq2.2} we have
\begin{eqnarray}
y'_1(x)^{(p-1)}-y'_2(x)^{(p-1)}&=&(p-1)\int_0^xq(t)[y_1(t)^{(p-1)}-y_2(t)^{(p-1)}]dt\nonumber\\
&+&(p-1)\int_0^xw(t)[f(y_2(t))-f(y_1(t))]dt.\label{eq5.6}
\end{eqnarray}
Apply an application of  the mean value theorem: for $a_1$ and $a_2$ of the same sign,
\begin{equation}
\label{eq5.7}
 a_1^{(p-1)}-a_2^{(p-1)}=(p-1)(a_1-a_2)|\bar{a}|^{p-2}\ ,
\end{equation}
 where $\bar{a}$ lies between $a_1$, $a_2$. Let $W(x)=y_1(x)-y_2(x)$.

By \eqref{eq5.5} and \eqref{eq5.7}, we obtain that
\begin{eqnarray*}
\int_0^x|q(t)[y_1(t)^{(p-1)}-y_2(t)^{(p-1)}]|dt&\leq&(p-1)\|q\|_x\int_0^x|y_1(t)-y_2(t)||\bar{b}_t|^{p-2}dt\\
&\leq&2(p-1)\|q\|_x\int_0^x|W(t)||t|^{p-2}dt\\
&\leq&2(p-1)\|q\|_x\|W\|_x\int_0^x|t|^{p-2}dt,
\end{eqnarray*}
where $\bar{b}_t$ lies between $y_1(t)$ and $y_2(t)$ and recall that $\|\cdot\|_x$ means the maximum of the given function in the interval from $0$ to $x$ and $p-2>-1$. It follows the locally
Lipschitz continuity of $f$, and then
$$\int_0^x|w(t)[f(y_2(t))-f(y_1(t))]|dt\leq A\|W\|_x$$ for some positive constant $A$. So from \eqref{eq5.6}, we have that $$(p-1)|\bar{a}_x|^{p-2}|W'(x)|\leq B\|W\|_x,$$ where $\bar{a}_x$ is close to $\eta_2$ and B is some positive constant. That is, $|W'(x)|\leq C\|W\|_x$ near $0$, for some positive $C$.
By Lemma A, we have $W=0$ in a neighborhood of zero.
\item[(iv)] Let $\eta_1\neq 0$ and $\eta_2=0$. If $1< p\leq 2$, the right hand side of \eqref{eq2.2} is
locally Lipschitz continuous. So the local uniqueness is obvious.
Let $p> 2$. Recall the phase and radius function of the
Pr\"{u}fer substitution satisfying \eqref{eq3.2} and
\eqref{eq3.3}.
 Note that the radius function $r(x)> 0$ is
uniquely defined, while the phase function $\theta(x)$ is unique
modulo $2\pi_p$. Since the right hand side of \eqref{eq3.2} is
Lipschitz  in $\theta$ for $(x,\theta)\in [0,c]\times
[\pi_p/2,\phi]$ for some angle $\phi$, $\theta(x)$ is unique near
$\pi_p/2$.
 Hence, the radius function $r(x)$ is also unique. This
implies that $y$ is unique locally.
\end{enumerate}
By Proposition \ref{th2.1}, the solution $y(x ; \eta_{1}, \eta_{2})$ of the IVP
\eqref{eq1.1} and \eqref{eq2.1} is unique on $[0,1]$.
Finally, a general theory of the continuous dependence of the
solutions on initial conditions \cite[Chap. V, Theorem 2.1]{H82}
implies that $y(x ;
\eta_{1},
\eta_{2})$ and $y'(x ; \eta_{1}, \eta_{2})$ are continuous in $(x;\eta_{1}, \eta_{2} )\in [0,1]\times
\mathbb{R}^{2}$ (see also \cite{nt08}).

\section*{Acknowledgments}
 \hskip0.25in
The first author is supported in part by National Science Council,
Taiwan under contract number NSC 99-2115-M-022-001. The author
Y.H. Cheng is supported by National Science Council, Taiwan under
contract numbers NSC 98-2115-M-007-008-MY3.

\end{document}